\magnification=\magstep1
\input epsf.tex
\overfullrule=0pt

% font specification
\font\titlefont=cmbx8 scaled \magstep3
% Define Gothic Font from AMS family
%\font\teneufm=eufm10 scaled \magstep1 
%\newfam\eufmfam
%\textfont\eufmfam=\teneufm
%\def\frak#1{\fam\eufmfan#1}
\def\frak{\cal}

% local definitions
\def\r{{\vec r}}
\def\h{{\vec h}}
\def\p{{\vec p}}
\def\a{{\vec a}}
\def\hatsl{{\widehat{sl}_2}}
\def\P{{\frak P}}
\def\R{{\frak R}}
\def\L{{\frak L}}
\def\eq{\eqalign}
\def\eps{{\epsilon}}
\def\Z{{{\rm Sp}}}
\def\qed{\hfill $\Young 6pt (1)$}
\def\frak{\cal}
\def\sc{\scriptstyle}

%       frame2.tex     by M. Hashimoto          6/27/95
%
%  Usage.... \Young<size of a box>(partition)    (see an example below)   
%            \young<size of a box>(partition)    (for subscript)        
%
%\aho       ... temporary macro 
%\b@xsize   ... absolute size of one box in the Young diagram, default=6pt.
%               default=2.5pt for \young
%\p@rtition ... a sequence of nonnegative integers, separated by commas
%\getnum    ... get a number of boxes in the current row and put it into
%               \r@. it removes this number and the following comma (if any).
%               it must be \global.
%\if@cont   ... If \p@rtition is exhausted by \getnum it gets to be 
%               false it must be \global.
%\r@        ... The real value of \r@ is 2*(the number of boxes).
%               see \getnum. it must be \global.
%\Young     ... The macro that creates a Young diagram
%\@@hline   ... Creates a horizontal line.
%\r@wofbox  ... Creates a row of boxes.
%\@nd       ... Equivalent to &. But this token list is not forbidden.
%\@nds      ... Token list consisting of \@nd's.
\catcode`\@11 \newdimen\b@xsize
\newdimen\@plsize
\newtoks\p@rtition
\newcount\r@
\newtoks\@nds
\newif\if@cont
\def\Young#1(#2){\def\aho{#1}%
\def\t@ko{#2}\ifx\t@ko\empty\emptyset\else
\ifx\aho\empty\b@xsize6pt\else\b@xsize\aho\fi
\@plsize.2\b@xsize
\global\p@rtition{#2,\aho}\global\@conttrue\expandafter
\getnum\the\p@rtition
\global\setbox1\hbox to\b@xsize{\hfil
\vrule width0ptheight\b@xsize depth0pt}\vcenter{\hbox{\vbox{\offinterlineskip
\tabskip0pt\def\@nd{&}\halign{\vrule height
\b@xsize##&&\copy1##&\vrule##\cr\@@hline\thisline@\cr}}}}%
\fi}%
\def\young#1(#2){\def\aho{#1}%
\def\t@ko{#2}\ifx\t@ko\empty\emptyset\else
\ifx\aho\empty\b@xsize2.5pt\else\b@xsize\aho\fi
\@plsize.3\b@xsize
\global\p@rtition{#2,\aho}\global\@conttrue\expandafter
\getnum\the\p@rtition
\global\setbox1\hbox to\b@xsize{\hfil
\vrule width0ptheight\b@xsize depth0pt}\hbox{\vtop{\offinterlineskip
\tabskip0pt\def\@nd{&}\halign{\vrule height\b@xsize##&&\copy1##&\vrule
##\cr\@@hline\thisline@\cr}}}\fi}%

\def\thisline@{\crcr\r@wofbox\cr\@@hline
\if@cont\let\noriP\tr@ns\else\let\noriP\relax
\fi
\noriP}%
\def\tr@ns{\expandafter\getnum\the\p@rtition\thisline@}%
\def\getnum#1,#2\aho{\def\aho{#2}\ifx\aho\empty\def\aho{0}\global
\@contfalse\fi
\global\r@#1\global\multiply\r@\tw@
\global\p@rtition{#2\aho}}%
\def\@@hline{\omit
\mscount\r@\ifnum\r@=\z@\hbox to0pt{\hrulefill\hskip-\@plsize}\else
\loop\ifnum\mscount>0
\sp@n\repeat\hrulefill\fi
}
\def\r@wofbox{\global\@nds{}\global\count\@ne\r@
\boke\the\@nds}%
\def\boke{\ifnum\count\@ne>0\global\@nds
\expandafter{\the\@nds\@nd}\global\advance\count\@ne -\@ne
\global\let\ahoaho\boke\else\global\let\ahoaho\relax\fi
\ahoaho}%
\catcode`\@\active

\def\qed{\hfill $\Young 6pt (1)$}

\vbox{\vskip 4cm}
\centerline{\titlefont 
Spectral Decomposition
of Path Space}

\centerline{\titlefont in Solvable Lattice Model }

\par

\vskip 2cm

\centerline{
Tomoyuki Arakawa, 
Tomoki Nakanishi,
Kazuyuki Oshima, 
  and Akihiro Tsuchiya}

\bigskip

\centerline{Department of
Mathematics}

\centerline{
 Nagoya University} 

\centerline{Chikusa-ku,
Nagoya 464-01, Japan }

\vskip 3cm

\centerline{ABSTRACT}
We give the $spectral$ $decomposition$
of the path space of the
$U_q(\hatsl)$ vertex model with respect to the local energy functions. 
The result suggests 
the hidden Yangian module structure on the
$\hatsl$ level $l$  integrable
modules, 
which is consistent
with the earlier work  [1]
in the  level  one case.
Also we prove
the fermionic character
formula of the $\hatsl$ level $l$  integrable
representations in consequence.

\vskip 2cm

\eject

\centerline{\bf \S 1.\enskip Introduction}
\medskip
In the last decade of investigation, various close relations
between the solvable lattice model and the conformal field
theory have been revealed (for example, [2-5]).
The aim of this article is to point out a new 
interesting relation between the spectrum in the
solvable lattice model and the hidden quantum symmetry in the 
 conformal field theory.

Consider the higher spin vertex model  associated with
$l+1$ irreducible representation of $U_q(\hatsl)$ ([6,7]).
It is well-known that the characters of
the $\hatsl$ or $U_q(\hatsl)$ level $l$ integrable representations ${\frak L}(k)$
 can be calculated
by using its $path$ $space$ ${\frak P}(k)$ ([2,8]).
The energy of a path $\vec p$ is  given
by the sum of a sequence of numbers
$h(\vec p)=(h_1(\vec p),h_2(\vec p),...)$ minus the ground state energy
which depends on  the corresponding boundary condition. 
Here $h_i(\vec p)$ is the $i$-th local energy determined from 
the $i+1$-th component of $\vec p$ and its nearest neighbors by the local energy 
function. We propose the fact that the local
energy functions  not only play a combinatorial role , but also 
can be regarded as the $q\rightarrow 0$ limit of  the 
local integrals of motion which commutes with the corner transfer matrix.

At $q=0$, the energy of a path $\vec p$ is
essentially the eigenvalue of the logarithm of the corner transfer
 matrix action on  the one dimensional 
configuration space $\sum_{\vec p\in {\frak P}(k)}{\bf C}{\vec p}$. 
Hence $\vec p$ itself  is  the ^^ ^^ eigenvector" of the corner transfer matrix,  and at the same time  
it is a simultaneous
^^ ^^ eigenvector" of the  mutually commuting infinitely many 
^^ ^^ local  operators" $h_i$ at $q=0$.

In this paper we  studied the $spectral$ $decomposition$
of the path space with respect to the local energy functions $h_i$.
That is, we decomposed the path space ${\frak P}(k)$ as
$${\frak P}(k)=\coprod{\frak P}(k)_{\vec h}\,,$$
where ${\frak P}(k)_{\vec h}$ denotes the ^^ ^^ eigenspace" with 
the spectrum $\vec h$.
 We found that the spectrum can be
parameterized by  the restricted 
paths and
the Young diagrams.
Moreover, it turns out that the spectrum has  a rich degeneracy structure.
Determining  the explicit form of the $sl_2$ character of
${\frak P}(k)_{\vec h}$, we see that it  equals to that of 
an irreducible Yangian  module.

In the level one case, the Yangian structure on the $\hatsl$ integrable
 modules has been described in [1].
Being consistent with their work, what our result suggests is that
even in the higher level cases, there is a canonical
Yangian structure on the integrable module ${\frak L}(k)$
such that
$${\frak L}(k)\simeq \bigoplus_{\vec h}
W_{\vec h}\, ,$$
where each $W_{\vec h}$ is an irreducible
finite dimensional Yangian module
with ${\rm ch}W_{\vec h}={\rm ch}{\frak P}_{\vec h}$.
(A part of
the Yangian structure is written  in [9].)

We believe that this correspondence between the spectral
decomposition of the solvable lattice model and the ^^ ^^ Yangian
multiplets" in the conformal field theory is not a  mere coincidence,
but will play the important role in investigating the hidden quantum
symmetry in more general conformal field theories.

In addition, as a byproduct, the spectral decomposition 
leads  to the fermionic
 character formula conjectured in
[9] and proved independently in [10].

\smallskip
This paper is organized as follows.
In section 2
we review the basic facts we use in this article
about the vertex model.
In  section 3 we parameterize the spectrum of the path space
with the restricted paths and
 the Young diagrams.
In section 4 
the degeneracy of the spectrum is calculated.
In section 5 we drive several character
formulas.
In section 6 we discuss the connection of our results in the vertex model
with the hidden Yangian symmetry in the conformal field theory.
\bigskip

\centerline{\bf\S 2.\enskip Local energy function and character formula }
\medskip
In this section we review basic facts about the vertex model.
 See  [2,6,7,8,11]  for the detail.
We consider the higher spin vertex model of $U_q(\hatsl)$.  
Throughout this paper we fix an integer
$l\in {\bf N}$.
Let $S=\{l,l-2,\dots -l\}$ and
let $V$ be the $l+1$ dimensional irreducible representation
of $U_q(sl_2)$ with the standard basis  $\{v_s|\,
{\rm wt}(v_s)=s,\,s\in S \}$.
We can  extend this to the homogeneous evaluation representation
$V(z)$  of $U_q(\hatsl)$ with a 
complex parameter $z\ne 0$ . 
The $ R$-matrix $\check R(w)$ is the $U_q(\hatsl)$-intertwiner 
$V(z_1)\otimes V(z_2)\rightarrow V(z_2)\otimes V(z_1)$ [6], where $w={z_1\over z_2}$. 
 We put
$$\check R(w)=\sum_{s_1,s_2,s'_1,s'_2\in S_l} 
            \check R(w)_{s_1 s_2}^{s'_1,s'_2}\, 
          E_{s'_1 s_1} \otimes E_{s'_2 s_2} \, ,\eqno(2.1)  $$
where $E_{s s'}$ is an element in ${\rm End}_{\bf C}(V(z))$
such that
$E_{s s'}\cdot v_{s''}=\delta_{s',s''} v_s$ 
for $s,s',s''\in S$.

The $R$-matrix above satisfies
several properties, especially
the diagonal nature at $q=0$;
$$\eq{\lim_{q\rightarrow 0}\check R(w)
&=\sum_{s_1,s_2\in S_l} w^{-H(s_1,s_2)}
E_{s_1 s_1}\otimes E_{s_2 s_2}}\eqno(2.2)$$
up to a scalar multiplication, where  $H$ is a function
$$ H: S \times S \longrightarrow
         \{0,1,\dots,l\} $$
such that
$$ H(s,s') = 
   \cases{ {1 \over 2} (s'+l)
,&if $s+s'\ge0$ ,
       \cr {1 \over 2} (l-s),
&if $s+s' < 0$.\cr} \eqno(2.3)$$ 
We call this function the 
{\it H function} 
or the {\it local energy function}. 
By setting
$$ H(s,s') = H_{s s'}, $$ 
we can write the above definition in a matrix form;
$$ ( H_{ss'} ) =\bordermatrix{
                     &\scriptstyle{  l  } &\scriptstyle {l-2} & \dots &  &\scriptstyle{-l+2}&\scriptstyle{-l} \cr
\scriptstyle{l}    &l                  & l-1		              & \dots &  & 1                 &0\cr
\scriptstyle{l-2}   & l                & l-1              & \dots &  &1                & 1\cr
\vdots               &\dots               &                  &       &  &                & \cr
                     &\dots               &                  &       &  &\vdots             &\vdots\cr
\scriptstyle{-l+2} & l                 & l-1                & \dots &  &l-1                & l-1\cr
\scriptstyle{-l }  & l                  & l                & \dots &  &l                & l}\eqno (2.4)$$

Each matrix element
$\check R(w)_{s_1,s_2}^{s'_1,s'_2} $ 
defines a Boltzmann weight of the configuration
 ($s_1$,$s_2$,$s'_1$,$s'_2$)
round a vertex in a planner square lattice,
where $s_1$, $s_2$,  $s'_1$, and $s'_2$ take value in $S$.
It is often written as
$$\raise 20pt\hbox{$\check R(w)_{s_1,s_2}^{s'_1,s'_2}=$}
\epsfbox{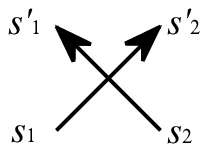}. $$
The corner transfer matrix method reduces the local state probability
to the one-dimensional configuration sum [7].
Fix a boundary condition $s_0\in S$ and let 
${\Sigma}_N=\{\vec s=(s_1,s_2,\cdots,s_N)| s_i\in S\}$.
The corner transfer matrix $A_N^{(s_0)}(w)$  is a matrix whose elements
are defined as
$$ \eq{A_N^{(s_0)}(w)_{\vec s}^{\vec s'}&
=\sum_{\rm configuration\enskip in\enskip the\enskip interior
       \enskip edges }\prod\check R(w)_{\eps_1,\eps_2}^{\eps'_1,\eps'_2}\cr
&={\atop \epsfbox{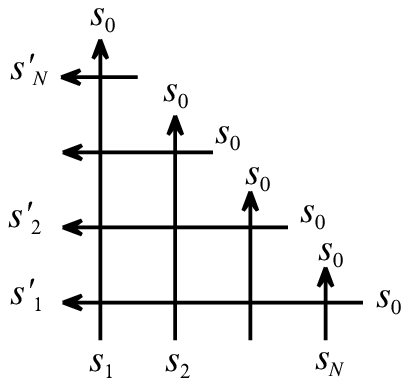}}}\eqno(2.5) $$
\vskip -0.5cm
\noindent for $\vec s=(s_1,\dots ,s_N), \vec s'=(s'_1,\dots ,s'_N) 
\in {\Sigma}_N$.

 Consider the limit $q\rightarrow 0$.
By the diagonal nature of the $R$-matrix, the corner transfer
matrix is also diagonalized at $q=0$, namely, 
$$A_N^{(s_0)}(w)_{\vec s}^{\vec s'}= \delta_{\vec s}^{\vec s'} 
w^{-(H(s_1,s_2)+ 2H(s_2,s_3)+\cdots+NH(s_{N},s_0))}\eqno(2.6) $$
with the normalization of the  $R$-matrix as in (2.2).
Hence
the value $H(s_1,s_2)+
   2H(s_2,s_3)+\dots+NH(s_{N},s_0)$ essentially contributes to
 the energy of a configuration $\vec s\in {\Sigma}_N$ with this boundary condition $s_0$. 

\medskip

Let us pass to the case when the lattice size $N$ is infinity.

For $k=0,\dots,l$, let $\vec s^{(k)}$ denote the $k$-th
ground state 
$(s_1^{(k)},s_2^{(k)},\dots\enskip  )=(l-2k,-(l-2k),l-2k,-(l-2k),\dots\enskip )$.
Let $${\Sigma}(k)=\{\vec s=(s_1,s_2,\dots\enskip  )|s_i\in S\enskip
{\rm for \enskip all\enskip} i,
                 \enskip  \vec s\approx \vec s^{(k)}\},
\eqno(2.7) $$
where $s\approx \vec s^{(k)}$ denotes the condition that $s_i=s_i^{(k)} $ except for
finitely many $i$'s.
We call ${\Sigma}=\bigsqcup_{k=0}^l{\Sigma}(k)$ the space
of the spin configurations of the vertex model associated to the $l+1$ dimensional
representation of $U_q(\hatsl)$.

Equivalently, we can define the configuration space above
using a path walking on the $sl_2$-weight lattice 
instead of a spin configuration.
Let 
$${\frak P}(k)=\{\vec p=(p_1,p_2,\dots\enskip )| p_i-p_{i+1}
\in S \enskip{\rm for \enskip all\enskip} i,\enskip
\vec p \approx \vec p^{(k)}\},\eqno(2.10) $$
where $\vec p^{(k)}=(p_1^{(k)},p_2^{(k)},\dots )=
(k,l-k,k,l-k,\dots\enskip )$ denotes the ground state path.
Then the one-to-one correspondence 
$$\eq{{\frak P}(k)&\simeq {\Sigma}(k) \cr
(p_1,p_2,\dots\enskip )&\leftrightarrow (p_2-p_1,p_3-p_2,\dots\enskip )}\eqno(2.11) $$
preserves the ground state for each $k$.
We call ${\frak P}=\bigsqcup_{k=0}^l {\frak P}(k)$ the path space of the vertex model
associated to the $l+1$ dimensional irreducible representation
of $U_q(\hatsl)$. 
We identify ${\Sigma}(k)$ with  ${\frak P}(k)$ hereafter
 \footnote{$^1$}{Sometimes a spin configuration $\vec s$ is also called a path in the literature 
\par  through
this correspondence.}.

Define a map
 $$\eq{h:\P(k)&\longrightarrow \{0,1,\dots,l\}^{\infty}\cr
  \p&\longmapsto h(\vec p)=(h_1(\p),h_2(\p),\dots\enskip)    , }$$  by setting 
 $$\eq{
h_i(\vec p)=H(p_{i+1}-p_i,p_{i+2}-p_{i+1})},\eqno (2.12)$$
where the function  $H$ is defined in (2.3).
The value $h_i(\p)$ is called  the $i$-th local energy of $\p\in\P(k)$.
Let $\h^{(k)}=h(\p^{(k)})=(k,l-k,k,l-k,\dots \enskip)$.

Now put the (total) energy of a path  $\vec p\in {\frak P}(k)$ as
$$\eq{E(\vec p)&=\sum_{i=1}^{\infty}i(h_i(\vec p)-h_i^{(k)} )  }\eqno(2.13) $$
(See (2.6)). Also define the $sl_2$ weight of a path $\vec p \in {\frak P}$ as
$${\rm Wt}(\vec p)=p_1.\eqno(2.14) $$
The following theorem
 is proved in [2]. (See also [8] for its relation with the crystal basis. )

\proclaim {Theorem 2.1}.
For $k=1,\cdots l$, let $\triangle(k)={k(k+2)\over 4(l+2)}$ and let
${\rm ch}\P(k)=q^{\triangle(k)}
\sum_{\vec p\in{\frak P}(k)}
q^{E(\vec p)}z^{{\rm Wt}(\vec p)}$.
Then,
$${\rm ch}\P(k)=
{\rm ch}_{{\frak L}(k)}(q,z) \eqno(2.15) $$
where ${\rm ch}_{{\frak L}(k)}(q,z)$
is the character of 
the $\hatsl$ or $U_q(\hatsl)$ level {\it l} integrable
module of the highest weight
$(l-k)\Lambda_0+k\Lambda_1$ .

\par
\medskip

The image of $\P(k)$ by the map $h$ is denoted by Sp$(k)$.
Then,
\proclaim{Proposition 2.2}.
For $k=0,1, \dots l$,
 $$ \eq{{\rm Sp}(k)=\{\h=(h_1,h_2,\dots,\enskip )|
(1).&\enskip({\rm Nearest\enskip neighbour\enskip condition})\quad  
 h_i+h_{i+1} \geq l\cr 
(2).&\enskip({\rm  Boundary\enskip condition})\quad
 \h\approx \vec h^{(k)}\}. }$$
\par
\hskip -0.7cm {\it Proof}.
~~Let $\p\in\P(k)$. Then $h(\p)$ trivially satisfies the condition (2) by definition.
Let us assume that $h_i(\p)=j$ for some $j\in\{0,\dots,l\}$.
Then by the definition of $H$ (formula (2.4)), $p_{i+2}-p_{i+1}\leq -l+2j$,
which means $h_{i+1}(\p)\geq l-j$. Hence $h(\P(k))\subset \Z(k)$.
The $h(\P(k))\supset \Z(k)$ part is similar to the oncoming proof  of Theorem 4.2:
Put $z=1$ in it.
\qed
\medskip
\proclaim {Remark 2.3}.
If $(h_{i-1}(\p),h_{i}(\p),h_{i+1}(\p))=(j,l-j,j)$ for 
some $j\in\{0,1,\dots,l\}$, then
$p_i=p_{i+1}+l-2j$ and $ p_{i+1}=p_{i+2}+2j-l$,
i.e., $p_{i}$ and $p_{i+1}$ are uniquely determined by $p_{i+2}$.\par
\medskip
We call  $ {\rm Sp}= \bigsqcup_{k=0}^{l}{\rm Sp}(k) $
 the spectrum of the path space
${\frak P }$ with
respect to the local energy functions $h_i$.

\proclaim Main Problem.
Decompose the path space ${\frak P}$ by its spectrum with
respect to the local energy functions $h_i$.
\par

We carry out the {\it spectral decomposition}
of the path space 
in the next sections
.
\bigskip

\centerline{\bf\S 3. \enskip Decoding map}
\medskip

In this section we fix $k\in\{0,\dots l\}$.
By the nearest neighbor condition in Proposition 2.2, we can divide each $\h\in\Z(k)$ into segments 
$$\eqalign{\vec h=
    &(g_1|g_2|\dots|g_m|g_{\infty}) \cr
    &g_i=(h_{\gamma_{i}},h_{\gamma_{i}+1},\dots,h_{\gamma_{i+1}-1}),\enskip
i=1,\cdots m
\cr
 &g_{\infty}=(h_{\gamma_{m+1}},h_{\gamma_{m+2}},\dots\enskip)\cr
    &\gamma_1=1<\gamma_2<\gamma_3<\cdots<\gamma_{m+1}=M+1 } \eqno (3.1)$$ 
so that in each segment the sum of two adjacent elements
$h_i,h_{i+1}$  is always $l$, but any sum of the adjacent two lying across different segments
is  greater than $l$.
We call those segment $g_i$'s  elementary blocks of $\h$.
Let $l(g)$ denotes the number of the components in an elementary block $g$
 when it is finite.
Each $g_i$ has the form as follows:
$$\eq{g_i&=\left\{
\matrix{([[h_{\gamma_{i}}]]^{l(g_i)\over 2})\hskip 2.5cm {\rm if}~l(g_i)~{\rm is ~ even}\cr   
(h_{\gamma_{i}},\,[[l-h_{\gamma_{i}}]]^{l(g_i)-1\over 2})\hskip 1cm {\rm if}~l(g_i)~{\rm is ~ odd}}\right. 
\hskip 1cm {\rm for}~i=1,\dots,m,
     }\eqno(3.2)$$
where 
 $ [[i]]^{a} $ denotes
$$ \overbrace{i,l-i,i,l-i,\dots,i,l-i.}^{a~ {\rm pairs}}\eqno(3.3)  $$
And let us write the last block $g_{\infty}$ as
$$g_{\infty}=(h_{\gamma_{m+1}},[[l-h_{\gamma_{m+1}}]]^{\infty})=\left\{\matrix{
(k,l-k,k,l-k,\dots\enskip )\hskip 1cm {\rm if}\enskip M\enskip 
{\rm is\enskip even}\hfill\cr
(l-k,k,l-k,k,\dots\enskip )\hskip 1cm {\rm if}\enskip M\enskip 
{\rm is\enskip odd}\hfill  }\right.\eqno(3.4)$$
for later convenience.
Then, supposing there are $J$ elementary blocks of odd length,
we  can now rewrite the given $\h\in\Z(k)$  in the following form:
$$\eq{\h=(
&[[m_{11}]]^{b_{11}},[[m_{12}]]^{b_{12}},\dots\dots,[[m_{1n_1}]]^{b_{1n_1}},l_1,[[l-l_{1}]]^{c_{1}},\cr
&[[m_{21}]]^{b_{21}},[[m_{22}]]^{b_{22}},\dots\dots,[[m_{2n_2}]]^{b_{2n_2}},l_2,[[l-l_{2}]]^{c_{2}},\cr
&[[m_{31}]]^{b_{31}},[[m_{32}]]^{b_{32}},\dots\dots,[[m_{3n_3}]]^{b_{3n_3}},l_3,[[l-l_{3}]]^{c_{3}},\cr
&\dots.\hskip 6.9cm \dots,\cr
&[[m_{J1}]]^{b_{J1}},[[m_{J2}]]^{b_{J2}},\dots\dots,[[m_{Jn_J}]]^{b_{Jn_J}},l_J,[[l-l_{J}]]^{c_{J}},\cr
&[[m_{J+1\,1}]]^{b_{J+1\,1}},[[m_{J+1\,2}]]^{b_{J+1\,2}},\dots ,
[[m_{J+1\,n_{J+1}}]]^{b_{J+1\,n_{J+1}}},l_{J+1},[[l-l_{J+1}]]^{\infty}),}
\eqno(3.5) $$
where $$\eq{\{m_{ij} \} &:~{\rm the ~initial~ elements~ of~ the~ blocks ~of~ even ~ length,}\cr
\{l_i|1\leq i\leq J\} &:~{\rm the ~initial~ elements~ of~ the~ blocks ~of~ odd ~ length,} \cr
l_{J+1}=h_{\gamma_{m+1}}&:~{\rm the ~initial~ elements~ of~ the~ last~block ~}g_{\infty},   }
\eqno(3.6)$$
 and $\{b_{ij}\}$(resp. $\{c_i\}$) are some positive (resp. non-negative)  integers determined by the lengths of 
the corresponding blocks. Notice that $J\equiv M$ mod 2.

\proclaim{Example 3.1}. $(l=3,k=1)$ 
$$\eq{(a)~  
\h&=(1,2,1,2|2,1|3|0,3|2,1,2|3,0,3|1,2|3,0,3|1,2,1,2,1,2,\dots,\enskip ) \cr   
&=[[1]]^2,[[2]]^1,3,[[0]]^1,2,[[1]]^1,3,[[0]]^1,[[1]]^1,3,[[0]]^1,1,[[2]]^{\infty}\cr
(b)~\h&=(0,3|1,2,1|3,0,3|1,2,1,2,1,2|2,1,2,1|3,0,3|2,1,2,1,2,1,2,1,\dots,\enskip ) \cr
&=[[0]]^1,1,[[2]]^1,3,[[0]]^1,[[1]]^3,[[2]]^2,3,[[0]]^1,2,[[1]]^{\infty} }$$
\par
\medskip

\proclaim {Lemma 3.2}.
For a given $\h\in\Z(k)$,
consider the sequence of the {\it initial} elements of elementary blocks
$$\eq{\big(m_{11},m_{12},\dots,m_{1m_1},l_1,m_{21},\dots,m_{2n_2},l_2,m_{31},
\dots,m_{3n_3},l_3,\dots \dots,\dots m_{J+1\,n_{J+1}},l_{J+1}  \big) }\eqno(3.7)$$ 
in the notation of (3.5) and (3.6). Then,\vskip 0.1cm
$$\eqalign{ &0\leq m_{11} < m_{12} < \cdots < m_{1n_1} <l_1
\cr &l-l_{i-1} < m_{i1} < m_{i2} < \cdots < m_{in_i} <l_i,\quad
     {\rm for \enskip  }i=2,\dots,J+1. } \eqno (3.8) $$
\par
\hskip -0.7cm {\it Proof}.~~ The statement directly follows from the fact that 
  the initials of the adjacent elementary blocks $g_i=(h_{\gamma_i},\dots, )$ and
$g_{i+1}=(h_{\gamma_{i+1},\dots)}$ satisfy
$$\eq{ h_{\gamma_i}< h_{\gamma_{i+1}}\hskip 1cm {\rm if}~l(g_i)~{\rm is ~ even} \cr
 l-h_{\gamma_i}<h_{\gamma_{i+1}}  \hskip 1cm {\rm if}~l(g_i)~{\rm is ~ odd}  }\eqno(3.9)$$
by Eq.(3.2) 
\qed
\medskip
For the given $\h\in\Z(k)$, let
 $$\h^{\sharp}=(l_1,\dots,l_J,l_{J+1})\eqno(3.10)$$ be the numbers defined in (3.6).
From the lemma above, 
we see that   a sequence of the form
$$ (l_i), [[l-l_i]]^{d_1},  [[l-l_i+1]]^{d_2}, [[l-l_i+2]]^{d_3}, \dots,
[[l_{i+1}-1]]^{d_{l_{i+1}-l+l_i}}, (l_{i+1})    \eqno(3.11)  $$
with $d_i\in {\bf Z}_{\geq 0}$
is inserted  between two initials $\l_i,l_{i+1}$ in  $\h$.
Set
$$\eq{  N &=N(\h)= \sum_{j=0}^J(l_{j+1}-l+l_j)\cr
t_i&=t_i(\h)=\sum_{j=0}^{i-1}(l_{j+1}-l+l_j)\hskip 1cm {\rm for}~ i=1,\dots J,}\eqno(3.12)$$
where $l_0=l$, and $t_1<t_2<\dots <t_J<N$ by the  above lemma.  
Then we can find a sequence
$$a(\h)=(a_1,\dots,a_N)\in({\bf Z}_{\geq 0})^N \eqno(3.13)$$
so that
$$\eq{\h=\big(&[[0]]^{a_1},[[1]]^{a_1},\dots,[[l_1-1]]^{a_{t_1}},l_1,\cr
&[[l-l_1]]^{a_{t_1+1}},[[l-l_1+1]]^{a_{t_1+2}},\dots,[[l_2-1]]^{a_{t_2}},l_2,\cr
& \dots\hskip 7cm \dots,\cr
&\hskip 5.5cm\dots,[[l_i-1]]^{a_{t_i}},l_i,\cr
&[[l-l_i]]^{a_{t_i+1}}[[l-l_i+1]]^{a_{t_i+2}}\dots,\cr 
&\dots\hskip 7cm\dots, \cr
&\hskip 5.5cm\dots,[[l_J-1]]^{a_{t_J}},l_J,\cr
& [[l-l_J]]^{a_{t_J+1}},\dots \hskip 3cm\dots,[[l_{J+1}-1]]^{a_N},l_{J+1},
[[l-l_{J+1}]]^{\infty}\big).  }\eqno(3.14)$$

We regard the above defined $a(\h)$ as  a  Young diagram  of depth $N=N(\h)$ by
$$\eq{ ({\bf Z}_{\geq 0})^N&\simeq Y_N=\{(\lambda_1,\lambda_2,\dots,\lambda_N)|
\lambda_i\in{\bf Z},\,\lambda_1\geq \lambda_2\geq\dots \lambda_N\geq 0 \} \cr
\a=(a_1,\dots,a_N)&\mapsto (\sum_{i=1}^N a_i, \sum_{i=1}^{N-1} a_i,\dots,a_1),  }\eqno(3.15)$$
where $Y_N$ denotes the set of the Young diagrams of depth $N$
(See Figure 1). For the later convenience, we set  $a(\h)=\phi$ when $N(\h)=0$,
and $Y_0=\{\emptyset\}$.
\topinsert$$\epsfbox{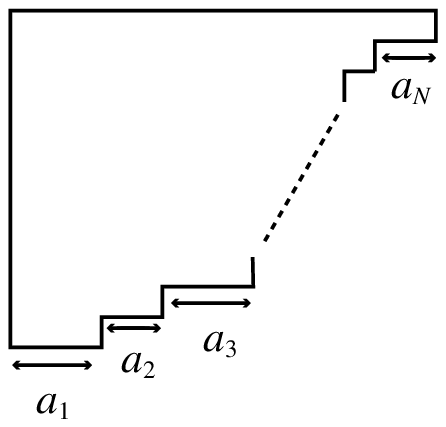}$$
\vskip -4.5cm$${\rm Figure~ 1.~~ Young ~ diagram~corresponding~to~}\a=(a_1,\dots,a_N)  $$
\endinsert
\medskip
\proclaim {Example 3.3}. (Continued from Example 3.1.)
$$\eq{(a)~ h&=[[0]]^0,[[1]]^2,[[2]]^1,3,[[0]]^1,[[1]]^0,2,[[1]]^1 ,[[2]]^0,3,[[0]]^1,[[1]]^1,[[2]]^0,
3,[[0]]^1,1,[[2]]^{\infty}}$$
 It reads as $N=11$ and 
$$\eq {a(\h)&=(0,2,1,1,0,1,0,1,1,0,1)\cr
& =\Young 5pt (8,7,7,6,5,5,4,4,3,2,0)}$$
$$\hskip -4.5cm\eq{(b)~\h&=[[0]]^1,1,[[2]]^1,3,[[0]]^1,[[1]]^3,[[2]]^2,3
,[[0]]^1,[[1]]^0,2,[[1]]^{\infty} }$$
It reads as $N=7$ and $$\eq{a(\h)&=(1,1,1,3,2,1,0)\cr
&=\Young 5pt (9,9,8,6,3,2,1)   }$$
\par
\medskip
We have seen so far that  once  $\h^{\sharp}$ is known,
then $\h\in\Z(k)$ is uniquely
determined by the index $ a(\h)\in Y_{N(\h)}$. 
Now let us consider the sequence $\h^{\sharp}=(l_1,\dots,l_J,l_{J+1})$ for $\h\in\Z(k)$.
We shall associate with it a restricted path of length $N=N(\h)$.

A sequence $\r=(r_0, \dots,r_N)$ is called a level $l$ restricted path of length $N$ if it 
satisfies the conditions
$$\eq{r_i\in\{0,\dots,l\}, ~
r_i-r_{i+1}=\pm 1.            } \eqno$$
Let $\R_N(k)$ denote the set of the level $l$ restricted paths of length $N$ with
the condition
$$r_0=0,~ r_N=k.  $$
Note that $\R_N(k)\ne \emptyset$ if and only if $N\equiv k$ mod 2 and $N\geq k$.
Since the sequence $\h^{\sharp}=(l_1,\dots, l_J,l_{J+1})$ satisfies 
$$0<l_1>l-l_2<l_3>\cdots <l_{2i-1}>l-l_{2i}<l_{2i+1}>\cdots \cases{>l-l_J<l_{J+1}=k~~
{\rm if}~J~{\rm is ~ even}\cr
<l_J>l-l_{J+1}=k~~{\rm if}~J~{\rm is ~ odd}} \eqno(3.16)$$ 
 by Lemma 3.2, we can define $r(\h)\in\R_{N}(k)$  for $\h\in\Z(k)$ as
$$\eq{&\scriptstyle
{ 0{\rm th}\hskip 1.9cm t_1{\rm  th} \hskip 3.2cm t_2 {\rm  th} 
}\cr
r(\h)=(&0,1,\dots,l_1-1,l_1,l_1-1,\dots,l-l_2+1,l-l_2,l-l_2+1,\dots\cr
 &\scriptstyle{\hskip 2.2cm  t_{2i-1} {\rm  th} \hskip 3.8cm t_{2i} {\rm  th} 
}\cr
&\dots,l_{2i-1}-1,l_{2i-1},l_{2i-1}-1,\dots,l-l_{2i}+1,l-l_{2i},l-l_{2i}+1,\dots,\cr
&\scriptstyle{\hskip 9cm N {\rm  th}
}\cr&\dots\hskip 8cm \dots,k).} \eqno(3.17) $$
(See Figure  2.)
Note that the numbers $(t_1,\dots,t_J)$ in (3.12) are interpreted as the extremal points of the
extremums $(l_1,l-l_2,l_3,\dots)$ of the path $r(\h)$.
\topinsert$$ \epsfbox{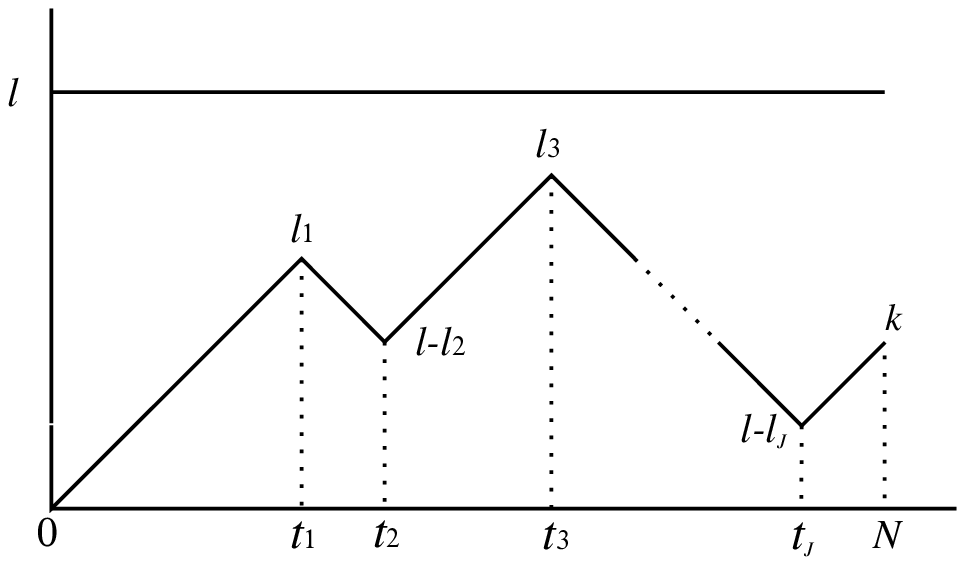}$$
\vskip -2cm $${\rm Figure~ 2. ~ ~restricted~ path }~r(\h) ~~(J:{\rm even}) $$
\endinsert
\medskip
\eject
\proclaim {Example 3.4}. (Continued from Example 3.3.)
\vskip 0.1cm
 (a) $\h^{\sharp}=(3,2,3,3,1).  $
So
$$\raise 50pt\hbox{$~~r(\h)=$} \hskip 0.5cm  \epsfbox{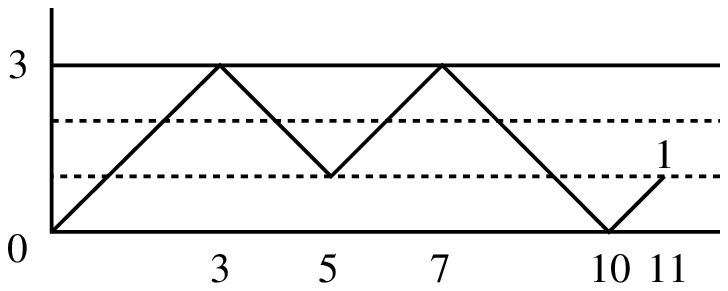} $$
(b) $\h^{\sharp}=(1,3,3,2).  $
So
$$ \hskip -2.2cm\raise 50pt\hbox{$~~r(\h)=$} \hskip 0.5cm \epsfbox{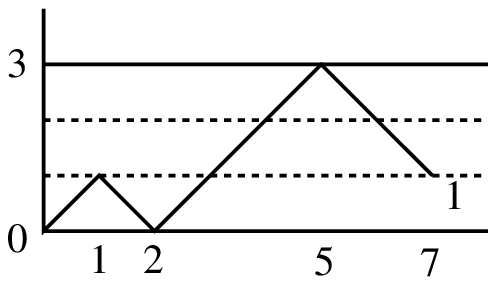}$$
\par
\medskip
Now let us summarize the above argument. We have found the map
$$\eq{  \pi_k : \Z(k)
 &\longrightarrow
 \bigsqcup_{N=k\atop N\equiv k~{\rm mod}~ 2}^{\infty}(\R_N(k) \times Y_N)\cr
\h&\longmapsto (r(\h),a(\h))}\eqno (3.18)$$
defined by (3.13) and (3.17).
We call the map $\pi_k$ the decoding map of $\Z(k)$.
\proclaim {Theorem 3.5}.
The decoding map
$\pi_k $
is bijective.
\par
\hskip -0.7cm {\it Proof}.~~ It is easy to construct the inverse map.
\qed
\medskip

\medskip
By the theorem above, 
we identify the spectrum Sp$(k)$ of the path space $\P(k)$
with $ {\displaystyle \bigsqcup_{N=k\atop N\equiv k~{\rm mod}~ 2}^{\infty}}$ 
$(\R_N(k) \times Y_N)$.
For $\vec r\in \R_N(k)$, $\vec a\in Y_N$,
let 
$$ \P(k)_{\r,\a}=\{\vec p\in\P(k) |~ h(\vec p)=\pi_k^{-1}(\vec r, \vec a) \}\eqno (3.19)$$
be the path space of the spectrum $(\r,\a)$. And let
$${\rm ch}\P(k)_{\vec r, \vec a}=
\sum\limits_{\vec p\in\P(k)_{\r,\a}}
q^{\triangle(k)+E(\p)}z^{{\rm Wt}(\vec p)}.\eqno (3.20)$$
be the character of $\P(k)_{\r,\a}$.
Then by definition
$$ {\rm ch}\P(k)=\sum_{(\r,\a)\in \Z(k)}{\rm ch}\P(k)_{\r,\a}.\eqno(3.21)$$

 \bigskip

\centerline{\bf \S4.  Degeneracy of the spectrum}
\medskip

In this section we calculate the character  
${\rm ch}{\frak P}(k)_{\vec r, \vec a}$,
which describes the degeneracy of the spectrum.

For $ \vec{r} \in {\frak R}_N(k), $  we consider the sequence
$$ \vec{n} = ( n_1,\dots,n_N) , \quad n_1 \leq \cdots \leq n_N \eqno(4.1)$$
obeying the condition
$$ n_1=0, \;\; n_i = \cases{ n_{i-1}+1,& if $ r_{i-2}=r_i < r_{i-1}$,\cr
                        n_{i-1}, & otherwise.\cr}\eqno(4.2) $$
We define the {\it degree} of the restricted path $ \vec{r} \in R_N(k) $,
denoted by $ d(\vec{r}) $, as 
$$ d(\vec{r}) = \sum_{i=1}^N n_i. \eqno(4.3)$$

For an arbitrary $ \vec{a} \in Y_N $, we define the {\it size} of the
Young diagram, denoted by $ |\vec{a}| $, as
$$ |\vec{a}| = \sum_{i=1}^{N} (N+1-i)a_i. \eqno(4.4)$$

This is in accordance with the usual definition of the size of the Young diagram through
the correspondence (3.15).

\proclaim{Proposition 4.1}.
{Fix a given 
 $ (\vec{r},\vec{a}) \in {\frak R}_N(k) \times Y_N, \; 
   ( N \geq k,\; N \equiv k \;{\rm mod}\; 2)$. 
 Then for any $ \vec{p} \in {\frak P}(k)_{\vec{r},\vec{a}} $ the equality
 $$ E(\vec{p}) = d(\vec{r}) + |\vec{a}| \eqno(4.5)$$
 holds. }
\par

\hskip-0.7cm{\it Proof.} \quad  The proof will be given  by the induction for 
$ a_1,\ldots,a_N $.

If $ \vec{a} =\vec{0} $, it is easy to show the both sides of (4.5) are
$$   \cases{{\displaystyle \sum_{i=1}^{J}il_i+{J \over 2}k- {J(J+2) \over 4}l},\;
                   & if  $J$ is even, \cr
            {\displaystyle \sum_{i=1}^{J}il_i-{{J+1}\over 2}k- 
                {(J+1)(J-1) \over 4}l},\; 
                   & if  $J$ is odd.\cr
                       }  \eqno(4.6)$$

Next suppose that the statement is true for some  $ \vec{a} \in Y_N $.
Let $ \vec{a'} \in Y_N $ be an element such that 
$ a'_i = a_i+ \delta_{i,I} $ for some $I \; (1\leq I \leq N)$.
It is sufficient to prove that for any $ \vec{p'}\in {\frak P}(k)_{\vec{r},\vec{a'}} $ ,
$$ E(\vec{p'}) = E(\vec{p}) + N+1-I, \quad 
   \vec{p} \in {\frak P}(k)_{\vec{r},\vec{a}}.\eqno(4.7)$$

Let $ l_1,\dots,l_J,t_1,\dots,t_J $ be the numbers defined in (3.6) and (3.12). 
Suppose
$$ t_{m-1} < I \leq t_m \eqno(4.8)$$
for some $ m, (1 \leq m \leq J+1)$, where $ t_0=0,t_{J+1}=N $. We set
$$ I=t_{m-1}+t,\enskip 1 \leq t \leq t_m-t_{m-1}. \eqno(4.9)$$
Then $ \pi_k^{-1}(\vec{r},\vec{a'}) $ is obtained by inserting a term
$ [[l-l_{m-1}+t-1]] $ between the $(m+2(a_1+\cdots+a_I)-1) $-th component of
$ \pi_k^{-1}(\vec{r},\vec{a}) $ and the
one right next to.
Thus 
$$ \eqalignno{E(\vec{p'}) - E(\vec{p})
    &= \{(m+2(a_1+\dots+a_I))(l-l_{m-1}+t-1)+
            (m+2(a_1+\dots+a_I+1))(l_{m-1}-t+1)\} \cr
    & \phantom{abcde}+ \{ 2l(a_{I+1}+\dots+a_N)+2(l_m+\dots+l_J) \cr
    & \phantom{abcderv}-(2(a_1+\dots+a_N)+J+1)l_{J+1}+
            (2(a_1+\dots+a_N)+J+2)(l-l_{J+1})\},\cr
   } $$
where the first term is the contribution from the inserted 
sequence $  [[l-l_{m-1}+t-1]] $  and the second term is the one from the shift of 
the sequence following the inserted sequence. Hence
$$ \eqalignno{E(\vec{p'}) - E(\vec{p})
    &=l_{m-1}+2(l_m+\dots+l_{J})+l_{J+1}-(J+2-m)l-t+1\cr
    &=N-t_{m-1}-t+1\cr 
    &=N-I+1 .
    } $$ \qed

Let us turn to the calculation of the $ sl_2 $ part of the character, 
$$ 
\sum_{\vec{p} \in {\frak P}(k)_{\vec{r},\vec{a}}}z^{{\rm Wt}(\vec{p})}. 
\eqno(4.10)$$  
We first introduce a partition of $N$ associated to 
each Young diagram of depth $N$.
Let $ B_N $ be the set of the ordered partitions of $N$, i.e.,
$$ B_N = \{ \vec{b}=(b_1,\cdots,b_s) | \; b_1,\cdots,b_s \in {\bf N}, \;
            1 \leq s \leq N, \; \sum_{i=1}^{s}b_i= N \}. \eqno(4.11)$$
We define the map $ \beta : Y_N \longrightarrow B_N $ as follows: 
For $ \vec{a} = (a_1,\ldots,a_N) \ne \vec{0} \in Y_N $, let 
$ 2 \leq \beta_1 < \cdots < \beta_{s-1} \leq N $ be the numbers such that
$$ a_{\beta_1},\ldots,a_{\beta_{s-1}} \eqno(4.12)$$
exhaust all the nonzero elements in $ (a_2,\dots,a_N) $. ( Notice that we excluded 
the element $a_1$.) Then we define 
$$ \beta : 
   \vec{a} \mapsto \vec{b} = (\beta_1-1,\beta_2-\beta_1,\ldots,N+1-\beta_{s-1}) 
   \eqno(4.13)$$
If $ \vec{a} = \vec{0} $, we set $ \beta(\vec{a}) = (N) \in B_N $.
Note that the map $\beta $ is a surjection.
This definition is more easily seen in Figure 3.
\topinsert 
$$\epsfbox{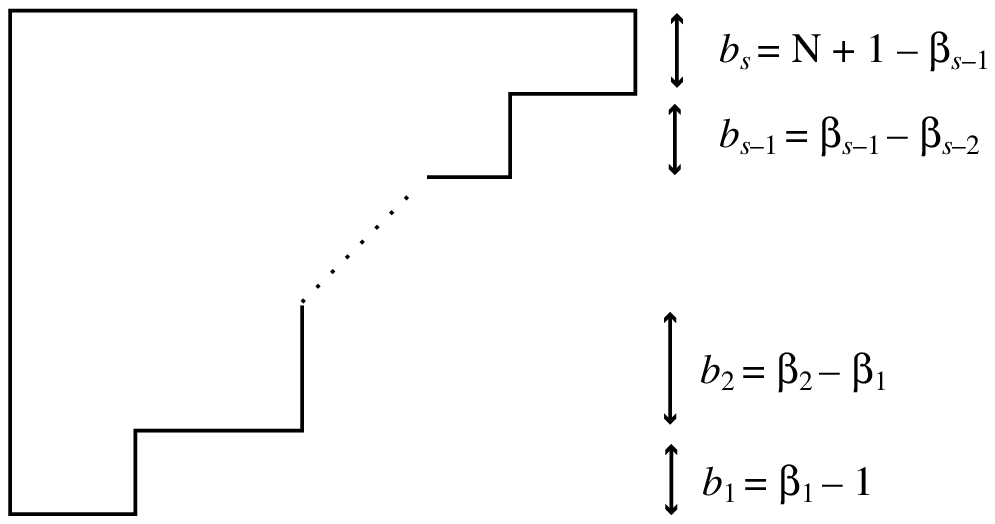} $$
\vskip-3.5cm
\centerline{ Figure 3.  The partition of $N$ associated to a Young diagram $\vec{a} \in Y_N $ }
\vskip-0.4cm
$$\hskip1.5cm  \beta(\vec{a}) = \vec{b} = (b_1,b_2,\dots,b_{s-1},b_s)$$
\endinsert

Let
$$ \chi_b(z) = {{z^{b+1}-z^{-b-1}} \over { z-z^{-1}}}. \eqno(4.14)$$
This is the character of the $(b+1)$-dimensional irreducible 
$ sl_2 $-module $ V_b$.
\proclaim{Theorem 4.2}.
{Let $ (\vec{r},\vec{a}) \in {\frak R}_N(k) \times Y_N, \;
       ( N \geq k, \; N \equiv k \; {\rm mod} \; 2)$.
 Let $ \vec{b} = (b_1,\ldots,b_s) \in B_N $ be the partition of $N$ associated
 to the Young diagram $ \vec{a} $ as in (4.13).
 Then 
$$ \sum_{\vec{p} \in {\frak P}(k)_{\vec{r},\vec{a}}}z^{{\rm Wt}(\vec{p})}
    = \prod_{i=1}^{s}\chi_{b_i}(z).\eqno(4.15)$$
 Combining (4.15) with (4.5), we have
$$ {\rm ch} {\frak P}(k)_{\vec{r},\vec{a}} =
                  q^{\Delta(k)+d(\vec{r})+|\vec{a}|} 
            \prod_{i=1}^{s} \chi_{b_i}(z). \eqno(4.16)$$
}
\par

The proof of Theorem 4.2 is divided into several steps.
 
(Step 1.) In this step we introduce the {\it incidence matrix } in order to 
calculate (4.10). 

For any $ (\vec{r},\vec{a}) \in {\frak R}_N(k) \times Y_N, \;
          (n \geq k,\; N \equiv k \;{\rm mod}\;2)$, let $ \vec{h} $ be the
corresponding element in the spectrum $ {\rm Sp}(k) $ (Theorem 3.5).
Let $M$ and $J$ be the numbers defined in (3.1) and (3.5) respectively, 
i.e., $J$ is the number of the elementary blocks of odd length in $\vec{h}$ and
$$ M = J + 2\sum_{i=1}^N a_i. \eqno(4.17)$$
Let $ \Sigma(k)_{\vec{r},\vec{a}}$ be the subset of $ \Sigma(k) $ corresponding
to $ {\frak P}(k)_{\vec{r},\vec{a}} $ through the  identity (2.11) .
For $ \vec{p} \in {\frak P}(k)_{\vec{r},\vec{a}} $, 
let $\vec{s} \in \Sigma(k)_{\vec{r},\vec{a}} $ be the spin configuration 
corresponding to $ \vec{p}$.    
Then it has the following form:
$$ \vec{s} 
 = (s_1,\dots,s_{M+1},l-2l_{J+1},-l+2l_{J+1},l-2l_{J+1},\dots). \eqno(4.18)$$
The summation $ s_1 + \cdots + s_{M+1} $ represents how the path goes up or down 
from the starting point $p_1$ to the $ (M+2) $-th point of $\vec{p}$, i.e.,
$l-l_{J+1} $. 
Therefore the following equality holds:
$$ \sum_{\vec{p} \in {\frak P}(k)_{\vec{r},\vec{a}}}z^{{\rm Wt}(\vec{p})}
    =  z^{l-l_{J+1}} \sum\limits_{\vec{s} \in \Sigma(k)_{\vec{r},\vec{a}}} 
              z^{-(s_1+\cdots+s_{M+1})}.  \eqno(4.19)$$
We define the function $ F(\vec{r},\vec{a};z) $ as
$$ F(\vec{r},\vec{a};z) = \sum_{\vec{s} \in \Sigma(k)_{\vec{r},\vec{a}}}
    z^{s_1 + \cdots + s_{M+1} }.\eqno(4.20)$$
This is evaluated by the multiplication of the following  $(l+1) \times (l+1)$ 
matrices: For $ a = 0,1,\cdots,l,$

$$ M_a =  \bordermatrix{
     &       &{\scriptstyle l-a+1}&          &        &           &\cr
     &       &     z^l &          &        &           &\cr
     &       & z^{l-2} &          &        &           &\cr
     &       & \vdots  &          &        &           &\cr
   {\scriptstyle a+1} &       & z^{l-2a}& z^{l-2a} & \cdots & z^{l-2a}  &\cr
                   &       &         &          &        &           &\cr },
     \eqno(4.21) $$

$$  H_a =\bordermatrix{
            &    &    & {\scriptstyle l-a+1} &          &        &          \cr
            &    &    &        &          &        &          \cr
{\scriptstyle a+1}&    &    &z^{l-2a}& z^{l-2a} & \cdots &  z^{l-2a}\cr
            &    &    &        &          &        &          \cr 
            &    &    &        &          &        &          \cr}, \;
   V_a =  \bordermatrix{
           &    &{\scriptstyle l-a+1}  &       &         &   &\cr
           &    &   z^l   &       &         &   &\cr
           &    & z^{l-2} &       &         &   &\cr
           &    & \vdots  &       &         &   &\cr
  {\scriptstyle a+1} &    & z^{l-2a}&       &         &   &\cr
           &    &         &       &         &   &\cr }. \eqno(4.22)$$
The matrix $ M_a $ is named the {\it incidence matrix }. 
This is defined in such a way that 

1) $ (M_a)_{ij} \ne 0 $ if and only if $ H_{l+2-2i,l+2-2j} \ne 0 $ in (2.4).

2) If $ (M_a)_{ij} \ne 0 $, then $ (M_a)_{ij} = z^{l+2-2i} $.

\hskip-0.7cm From these properties, the $ (i,j) $-component of the matrix
$ M_{h_1} \cdots M_{h_m} $
equals to 
$$ \sum_{\vec{s}=(s_1,\dots,s_m)} z^{s_1 + \cdots + s_m}, \eqno(4.23)$$
where the summation is taken over all the spin configurations $ \vec{s} $ such that
$$ H(s_i,s_{i+1}) = h_i, \; s_1=l+2-2i, \; s_{m+1}=l+2-2j. \eqno(4.24)$$

\proclaim{Lemma 4.3}.
{ 1). For an arbitrary $ \vec{h}={\pi_k}^{-1}(\vec{r},\vec{a}) $,
 the function $ F(\vec{r},\vec{a};z) $ is given by the sum of the matrix 
elements on the $ l-l_{J+1}+1 $-th column of the matrix
$$ M_{h_1} \cdots M_{h_M}M_{l_{J+1}}. \eqno(4.25)$$ \break
  2). Equivalently, the function $ F(\vec{r},\vec{a};z) $ 
is also given by the $ (l+1,l-l_{J+1}+1) $ component of the matrix
$$ T(\vec{r},\vec{a}) = z^l H_l M_{h_1} \cdots M_{h_M} V_{l_{J+1}}. \eqno(4.26)$$
}
\par
\hskip-0.7cm{\it Proof.} \quad
1). For $ \vec{s} \in \Sigma(k)_{\vec{r},\vec{a}}, \; s_{M+2}=-l+2l_{J+1} $ 
as in (4.18).
Applying the fact above Lemma 4.3, we have the statement.

\hskip-0.7cm 2). In (4.25) one can replace $M_{l_{J+1}}$ by $V_{l_{J+1}}$  
without changing the $l-l_{J+1}+1$-th column, which is the only part we need. 
It is clear from the explicit form of $ z^l H_l $ that the
$ (l+1,l-l_{J+1}+1) $ component of the matrix $T(\vec{r},\vec{a})$ is equal to
the sum of the matrix elements on the $l-l_{J+1}+1$-th column of (4.25).   \qed

\hskip-0.7cm{\bf Example 4.4.} \enskip
Here we consider the case
$$ \eqalign{ &l=3,\; N=3,\; k=1,\cr
            & \vec{r}=(0,1,2,1),\cr
            & \vec{a}=(0,1,0).} $$
In this case $J=1$, $M=3$, $l_{J+1}=2$, and $(l+1,l-l_{J+1}+1)=(4,2)$.
From these data, $\vec{h} = {\pi_k}^{-1}(\vec{r},\vec{a})$ is
$$ \vec{h} = (1,2,2,2,[[1]]^{\infty}). $$
The spin configurations in $ \Sigma(k)_{\vec{r},\vec{a}} $ are 
$$ (s_1,s_2,s_3,s_4,1,-1,1,-1,1,\ldots), $$
where
$$  \eqalignno{(s_1,s_2,s_3,s_4,) 
             &= (3,-1,1,1,),\cr
             &\phantom{=1} (3,-1,-1,1),\cr
             &\phantom{=1} (3,-1,-1,-1),\cr
             &\phantom{=1} (1,-1,1,1),\cr
             &\phantom{=1} (1,-1,-1,1),\cr
             &\phantom{=1} (1,-1,-1,-1). } $$
Thus
$$ F(\vec{r},\vec{a};z) = z^4 + 2 z^2 + 2 + z^{-2}. $$
On the other hand,
$$ \eqalignno{ T(\vec{r},\vec{a}) &= z^3 H_3 M_1 M_2 M_2 V_2  \cr
    & = z^3 \pmatrix{ 0 & 0 & 0   & 0 \cr
                  0 & 0 & 0   & 0 \cr
                  0 & 0 & 0   & 0 \cr
             z^{-3} & z^{-3} & z^{-3}   & z^{-3} \cr}
        \pmatrix{ 0 & 0 & z^3 & 0 \cr
                  0 & 0 & z   & z \cr
                  0 & 0 & 0   & 0 \cr
                  0 & 0 & 0   & 0 \cr}
        \pmatrix{ 0 & z^3 & 0 & 0 \cr
                  0 & z   & 0 & 0 \cr
                  0 & z^{-1} & z^{-1} & z^{-1} \cr
                  0 & 0 & 0 & 0 \cr} \cr
 & \phantom{Hello. How are you?} \times    \pmatrix{ 0 & z^3 & 0 & 0 \cr
                  0 & z   & 0 & 0 \cr
                  0 & z^{-1} & z^{-1} & z^{-1} \cr
                  0 & 0 & 0 & 0 \cr}
        \pmatrix{ 0 & z^3 & 0 & 0 \cr
                  0 & z   & 0 & 0 \cr
                  0 & z^{-1} & 0 & 0 \cr
                  0 & 0 & 0 & 0 \cr} \cr
    & = \pmatrix{ 0 & 0 & 0 & 0 \cr
                  0 & 0 & 0 & 0 \cr
                  0 & 0 & 0 & 0 \cr
                  0 & z^4+2z^2+2+z^{-2} & 0 & 0 \cr}. }$$
We see that the $(4,2)$ component of 
the matrix $ T(\vec{r},\vec{a}) $ is equal to the $ F(\vec{r},\vec{a};z) $. 

\bigskip
\medskip
 
Let us proceed to the evaluation of $ T(\vec{r},\vec{a}) $ in Lemma 4.3.
\smallskip
(Step 2.) In this step, we represent $ T(\vec{r},\vec{a}) $ as a factorized form (4.36). 

The next lemma is easily proved and useful for our purpose.

\proclaim{Lemma 4.5}.
{ For $ a= 0,1,\cdots,l$ and $ n \in {\bf N} $,
$$ (M_a M_{l-a})^{n} = V_a H_{l-a}. \eqno(4.27)$$
}
\par

\proclaim{Lemma 4.6}.
 { For $ \vec{a}, \vec{a'} \in Y_N $, if 
   $ \beta(\vec{a}) = \beta(\vec{a'}) $, then
$$ F(\vec{r},\vec{a};z) = F(\vec{r},\vec{a'};z). \eqno(4.28)$$
 }
 \par
\hskip-0.7cm {\it Proof.} \quad
If $ a_1 \ne 0 $,  there is a factor 
$ z^l H_l (M_0 M_l)^{a_1} $ at the left end of $ T(\vec{r},\vec{a}) $
(see (3.14) and (4.27)). Since
$$ M_0 M_l = \pmatrix{
                 1 & 1 & \cdots & 1\cr
                   &   &        &  \cr
                   &   &        &  \cr
                   &   &        &  \cr}, \eqno(4.29)$$
$ z^l H_l (M_0 M_l)^{a_1} = z^l H_l .$ 
So it does not affect the function
$ F(\vec{r},\vec{a};z) $ whether $a_1$ is equal to zero or not.
Next  we suppose that
there is some $ i \; (2 \leq i \leq N) $ such that 
$ a_i \geq 1 , \; {a'}_i \geq 1, $ and $ a_i \ne {a'}_i $.
But from Lemma 4.5,
$$ (M_{a}M_{l-a})^{a_i} = (M_{a}M_{l-a})^{{a'}_i} = V_{a}H_{l-a}. \eqno(4.30)$$
Thus we have the result. \qed

From this lemma we could assume that 
$$ a_1=0, \enskip a_2,\cdots,a_N \in \{0,1 \} \eqno(4.31)$$
in order to prove Theorem 4.2 without losing generality.
Until the end of this proof we set this assumption.

For $ \vec{a} \ne \vec{0} $, 
let $ 2 \leq \beta_1 \leq \cdots \leq \beta_{s-1} <N $
be the numbers defined in (4.12), and $0 < t_1 < \cdots <t_J $ be the 
extremal points of $ \vec{r} $. 
We align $ \beta_i$'s and $ t_i $'s as follows:
$$ \eqalign{
    & 0 < t_1 < \cdots < t_{m_1} < \beta_1-1 \leq t_{m_1+1} < \cdots <
      t_{m_2} < \beta_2-1 \leq \cr
    & \cdots \cdots \leq t_{m_{s-2}+1} < \cdots < t_{m_{s-1}} <
      \beta_{s-1}-1 \leq t_{m_{s-1}+1} < \cdots < t_J < N. \cr
      } \eqno(4.32)$$
Accordingly $ T(\vec{r},\vec{a}) $ has a form
$$ \eqalign{
    & z^l H_l M_{l_1} \cdots M_{l_{m_1}} (M_{h'_1}M_{l-h'_1})
      M_{l_{m_1+1}} \cdots M_{l_{m_2}} (M_{h'_2}M_{l-h'_2})
      M_{l_{m_2+1}} \cdots \cr
    & M_{l_{m_{s-2}}} (M_{h'_{s-2}}M_{l-h'_{s-2}}) M_{l_{m_{s-2}+1}} \cdots
      M_{l_{m_{s-1}}} (M_{h'_{s-1}}M_{l-h'_{s-1}}) M_{l_{m_{s-1}+1}} \cdots 
      M_{l_J} V_{l_{J+1}} , \cr }\eqno(4.33) $$
where 
$$ h'_i = \cases{ r_{\beta_i-1}   & if $m_i$ is even,\cr
                  l-r_{\beta_i-1} & if $m_i$ is odd. } \eqno(4.34)$$

Using Lemma 4.5 one can replace the factors $ (M_{a}M_{l-a}) $ by
$(V_{a}H_{l-a})$. Thus we have
$$ \eqalignno{T(\vec{r},\vec{a}) 
   &=z^l (H_l M_{l_1} \cdots M_{l_{m_1}}V_{h'_1})
     (H_{l-h'_1} M_{l_{m_1+1}} \cdots M_{l_{m_2}} V_{h'_2}) \cdots \cr
      &\phantom{=S}(H_{l-h'_{s-2}}M_{l_{m_{s-2}+1}} \cdots M_{l_{m_{s-1}}} V_{h'_{s-1}})
     (H_{l-h'_{s-1}}M_{l_{m_{s-1}+1}} \cdots  M_{l_J} V_{l_{J+1}}). &(4.35)\cr } $$
Therefore $ T(\vec{r},\vec{a})$ is a product of the factors
$$ \eqalign{ T(\vec{r},\vec{a}) &= z^l S_1 S_2 \cdots S_s ,\cr
  S_i &= H_{l-h'_{i-1}}M_{l_{m_{i-1}+1}} \cdots M_{l_{m_i}}V_{h'_i},} \eqno(4.36)$$
where
$$  m_0 = 0 , \;\;m_s = J, \quad \beta_0 =1, \quad \beta_s=N+1 . $$
\hskip-0.7cm {\bf Example 4.7.} \enskip
To illustrate the procedures (4.32 - 35), let us give one more example. Let
$$ \eqalign{
   & l = 3, \; N=9, \; k=1, \cr
   & \vec{r} = (0,1,2,1,2,3,2,1,0,1) \in {\frak R}_9(1), \cr
   & \vec{a} = (0,1,0,1,0,0,1,0,0) \in A_9. } $$
The sequence of the local energies is
$$ {\pi_1}^{-1}(\vec{r},\vec{a}) = ([[1]],2,2,[[1]],3,[[1]],3,[[1]]^{\infty}). $$
The Young diagram corresponding to $ \vec{a} $ is
$$ \Young 7pt (3,3,3,2,2,2,1,1,0). $$
Then we see that
$$ \eqalign{ \beta_1 = 2, \; \beta_2 = 4, \; \beta_3 = 7 \cr
             t_1 = 2, \; t_2 = 3, \; t_3 = 5, \; t_4 = 8,} $$
and the alignment of (4.31) in this case is
$$0 <1 <2 <3 \leq 3 <5 <6 <8. $$
Accordingly 
$$T(\vec{r},\vec{a};z)= z^3 H_3M_1M_2M_2M_2M_1M_2M_3M_1M_2M_3V_1. $$
From Lemma 4.5 this can be factorized as 
$$ z^3 (H_3V_1)(H_2M_2M_2V_1)(H_2M_3V_1)(H_2M_3V_1). $$ 

\bigskip

(Step 3.) Finally we evaluate the factor $ S_i $ explicitly. 
  
Let $ E_{i,j} $ be the matrix whose matrix element is $1$ at the $ (i,j) $ component 
and zero otherwise.

\proclaim{Lemma 4.8}.
{ $$ S_i = z^{h'_{i-1}-h'_i} \; \chi _{\beta_i - \beta_{i-1}}(z) \;
             E_{l-h'_{i-1}+1,l-h'_i+1} \eqno(4.37)$$ 
}
\par

\hskip-0.7cm{\it Proof.} \quad   
Let us put $ m_{i-1}+1=n,\enskip m_i =m $, for simplicity. We need to calculate
the product
$$  H_{l-h'_{i-1}}M_{l_n} \cdots M_{l_m} V_{h'_i}. \eqno(4.38)$$
First we consider the product $  H_{l-h'_{i-1}}M_{l_n} $.
$$ \eqalignno{ H_{l-h'_{i-1}}M_{l_n}
  &=  \bordermatrix{
                  &    &{\sc h'_{i-1}+1}    &       &                 \cr
                  &    &                 &       &                  \cr
{\sc l-h'_{i-1}+1}&    & z^{-l+2h'_{i-1}}& \cdots& z^{-l+2h'_{i-1}} \cr
                  &    &                 &       &                  \cr
                  &    &                 &       &                  \cr
                  &    &                 &       &                  \cr }
     \bordermatrix{
            &    &{\sc l-l_n+1}&        &            \cr
            &    & z^l         &        &            \cr
            &    & z^{l-2}     &        &            \cr
            &    & \vdots      &        &            \cr
 {\sc l_n+1}&    & z^{l-2l_n}  & \cdots & z^{l-2l_n} \cr
            &    &             &        &            \cr } \cr
 &=  \bordermatrix{
                  &    &{\sc l-l_n+1}            &   &        &   \cr
                  &    &                         &   &        &   \cr
{\sc l-h'_{i-1}+1}&    &\varphi_{l_n-h'_{i-1}}(z)& 1 & \cdots & 1 \cr
                  &    &                         &   &        &   \cr
                  &    &                         &   &        &   \cr }
                    z^{2(h'_{i-1}-l_n)}, &(4.39)\cr } $$
where 
$$ \varphi_m(z) = z^{2m} +z^{2m-2} + \cdots + 1. \eqno(4.40)$$
Since 
$ l_n-h'_{i-1} = t_n - (\beta_{i-1} - 1),$
$$  H_{l-h'_{i-1}}M_{l_n}
   =   \bordermatrix{
                  &    &{\sc l-l_n+1}            &   &        &   \cr
                  &    &                         &   &        &   \cr
{\sc l-h'_{i-1}+1}&    &\varphi_{t_n-(\beta_{i-1}-1)}(z)& 1 & \cdots & 1 \cr
                  &    &                         &   &        &   \cr
                  &    &                         &   &        &   \cr }\;
                    z^{-t_n+(\beta_{i-1}-1)+h'_{i-1}-l_n}. \eqno(4.41)$$
In a similar way
$$ \eqalign{ 
    &H_{l-h'_{i-1}}M_{l_n}M_{l_{n+1}} \cr
    &=  \bordermatrix{
                  &    &{\sc l-l_n+1}            &   &        &   \cr
                  &    &                         &   &        &   \cr
{\sc l-h'_{i-1}+1}&    &\varphi_{t_n-(\beta_{i-1}-1)}(z)& 1 & \cdots & 1 \cr
                  &    &                         &   &        &   \cr
                  &    &                         &   &        &   \cr }
        \bordermatrix{
            &    &{\sc l-l_{n+1}+1}&        &            \cr
            &    & z^l         &        &            \cr
            &    & z^{l-2}     &        &            \cr
            &    & \vdots      &        &            \cr
 {\sc l_{n+1}+1}&    & z^{l-2l_{n+1}}  & \cdots & z^{l-2l_{n+1}} \cr
            &    &             &        &            \cr }\cr
   & \hskip1cm \times z^{-t_n+(\beta_{i-1}-1)+h'_{i-1}-l-n} \cr
   &=   \bordermatrix{
                  &    &{\sc l-l_{n+1}+1}            &   &        &   \cr
                  &    &                         &   &        &   \cr
{\sc l-h'_{i-1}+1}&    &\varphi_{t_n-(\beta_{i-1}-1)+l_{n+1}-(l-l_n)}(z)& 1 & \cdots & 1 \cr
                  &    &                         &   &        &   \cr
                  &    &                         &   &        &   \cr }\;
              z^{-t_n+(\beta_{i-1}-1)+h'_{i-1}-l_n +l-2l_{n+1}}. \cr} \eqno(4.42) $$
Again, since 
$ l_{n+1}-(l-l_n) = t_{n+1}-t_n, $
$$ \eqalign{ 
    &H_{l-h'_{i-1}}M_{l_n}M_{l_{n+1}} \cr 
    &= \bordermatrix{
                  &    &{\sc l-l_{n+1}+1}            &   &        &   \cr
                  &    &                         &   &        &   \cr
{\sc l-h'_{i-1}+1}&    &\varphi_{t_{n+1}-(\beta_{i-1}-1)}(z)& 1 & \cdots & 1 \cr
                  &    &                         &   &        &   \cr
                  &    &                         &   &        &   \cr }\;
                    z^{-t_{n+1}+(\beta_{i-1}-1)+h'_{i-1}-l_{n+1}}. \cr } \eqno(4.43)$$
Repeating this procedure, the result is shown as 
$$ \eqalignno{ H_{l-h'_i}M_{l_n} \cdots M_{l_m} V_{h'_{i-1}}
   &=  \bordermatrix{
       &    &        &       &               {\sc l-h'_i+1}&   &\cr
       &    &        &       &                                &   &\cr
{\sc l- h'_{i-1}+1}&    &        &       &\varphi_{\beta_i-\beta_{i-1}}(z)&   &\cr
       &    &        &       &                                &   &\cr
       &    &        &       &                                &   &\cr} 
       z^{-\beta_i+\beta_{i-1}+h'_i-h'_{i-1}} \cr
   &=  z^{h'_i-h'_{i-1}} {\rm ch}_{V_{\beta_i - \beta_{i-1}}}(z)
       E_{l-h'_{i-1}+1, l-h'_i+1}.    &(4.44)}$$
\qed

\bigskip

Now we complete the proof of Theorem 4.2. 
By (4.20) we see
$$ \sum_{\vec{p} \in {\frak P}(k)_{\vec{r},\vec{a}}}z^{{\rm Wt}(\vec{p})}
    = z^{l-l_{J+1}} F(\vec{r},\vec{a};z^{-1}) \eqno(4.45) $$
From Lemma 4.8, $F(\vec{r},\vec{a};z)$ is evaluated as 
$$ F(\vec{r},\vec{a};z) =
   z^{l-l_{J+1}} \prod\limits_{i=1}^s \chi_{\beta_i - \beta_{i-1}}(z),
         \eqno(4.46)$$
where we have used 
$$ h'_0 - h'_s = -l_{J+1}. $$
Taking into the account of the definition of $ \beta(\vec{a}) $ ,
we get the formula
$$ \sum_{\vec{p} \in {\frak P}(k)_{\vec{r},\vec{a}}}z^{{\rm Wt}(\vec{p})}
    = \prod_{i=1}^s \chi_{b_i}(z),
    \quad \beta(\vec{a}) = ( b_1,\cdots,b_s). \eqno(4.47)$$ \qed

\bigskip

\centerline{ \bf \S5. \enskip Character formulas}
\medskip
In this section, as an application of our decomposition,  we derive 
some character formulas for the level $l$ integrable $ \widehat{sl}_2 $-modules.

Owing to Theorem 2.1, 3.5, 4.2 and  Propositions 4.1,  we immediately obtain 
the following formula :
\proclaim{Proposition 5.1}.
{$$ {\rm ch}_{{\frak L}(k)}(q,z) = q^{\Delta(k)} 
    \sum_{{N=k} \atop  N \equiv k \;{\rm mod} \; 2}^{\infty}    
    \sum_{(\vec{r},\vec{a}) \in {\frak R}_N(k) \times Y_N }
    q^{d(\vec{r})+|\vec{a}|} \prod_{i=1}^s \chi_{b_i}(z) ,
     \eqno(5.1)$$
where $(b_1,\cdots,b_s) = \beta(\vec{a})$. }
\par

Since the summation over $ {\frak R}_N(k) $ and $ Y_N $ are independent of 
each other, (5.1) can be rewritten as the following factorized form :

$$ \eqalign{
  {\rm ch}_{{\frak L}(k)}(q,z) &= q^{\Delta(k)} 
    \sum_{{N=k} \atop {N \equiv k \;{\rm mod} \; 2}}^{\infty}
                            F_{N,k}(q)G_N(q,z), \cr
  F_{N,k}(q) &= \sum_{\vec{r} \in {\frak R}_N(k)} q^{d(\vec{r})}, \cr
  G_N(q,z) &= \sum_{\vec{a} \in Y_N} q^{|\vec{a}|}
                 \prod_{i=1}^{s} \chi_{b_i}(z). } \eqno(5.2)$$

It is possible to evaluate the functions $ F_{N,k}(q) $ and 
$ G_N(q,z) $ more explicitly.
Let
$$ \eqalign{
    (q)_n &= \prod_{i=1}^n (1-q^i), \cr
    {N \brack n} &={ {(q)_N} \over {(q)_n (q)_{N-n}}}, 
                    \enskip {\rm for} \enskip 0 \leq n \leq N.} \eqno(5.3)$$

\proclaim{Proposition 5.2}.
{ $$ G_N(q,z) = {1 \over (q)_N} \sum_{n=0}^N {N \brack n}z^{N-2n}. \eqno(5.4)$$}
\par
\hskip-0.7cm{\it Proof.} 
Let $ \beta_i $'s be the numbers defined in (4.12). Then
$$ \eqalign{ G_N(q,z) 
     & = \sum_{\vec{a} \in Y_N} q^{|\vec{a}|}
                 \prod_{i=1}^{s} \chi_{b_i}(z) \cr
     & = \sum_{b_i \geq 1 \atop b_1+\cdots+b_s=N}
         \sum_{a_1 \geq 0 \atop a_{\beta_m} > 0}
          q^{\sum_{i=1}^{s}\{\sum_{m=0}^{i-1}a_{\beta_m}\}b_i}
             \prod_{i=1}^{i} \chi_{b_i}(z) \cr
    }\eqno(5.5) $$
Summed up with respect to $ a_{\beta_m} $'s, we see that the above equals to
$$ \sum_{b_i \geq 1 \atop b_1+\cdots+b_s=N}
   \prod_{i=1}^{s} \chi_{b_i}(z) 
   {q^{b_1} \over {1-q^{b_1}}}{q^{b_1+b_2} \over {1-q^{b_1+b_2}}}
   \cdots  {q^{b_1+\cdots+b_{s-1}} \over {1-q^{b_1+\cdots+b_{s-1}}}}
    {1 \over {1-q^{b_1+\cdots+b_{s}}}}.\eqno(5.6)$$
Considering it particularly as a summation  with respect to $ b_s$, 
again the above equals to
$$ \sum_{b_s=1}^{N} {q^{N-b_s} \over {1-q^N}} \chi_{b_s}(z)
     \sum_{b_i \geq 1 \atop b_1+\cdots+b_{s-1}=N-b_s}
     \prod_{i=1}^{s-1} \chi_{b_i}(z)
     {q^{b_1} \over {1-q^{b_1}}}{q^{b_1+b_2} \over {1-q^{b_1+b_2}}}
   \cdots  {1 \over {1-q^{b_1+\cdots+b_{s-1}}}}. \eqno(5.7)$$
Then we find that $ G_N(q,z) $ satisfies the following recursion relation:
$$ \eqalign{
    G_0(q,z) &= 1, \cr
    G_N(q,z) &= \sum_{b=1}^{N} {{q^b} \over {1-q^N}} 
                 \chi_{b}(z)G_{N-b}(q,z)   \cr } \eqno(5.8)$$
On the other hand, let 
$$ H_N(q,z) = \sum_{n=0}^N {N \brack n}z^{N-2n}.\eqno(5.9)$$
The function $ H_N(q,z) $ is essentially the Rogers-Szeg\"o
polynomial  [12], and satisfies the following recursion relation
$$ \eqalign{
    H_0(q,z) &= 1 ,\cr
    H_1(q,z) &= z + z^{-1} ,\cr
    H_N(q,z) &= (z+z^{-1}) H_{N-1}(q,z) - (1-q^{N-1})H_{N-2}(q,z).} \eqno(5.10)$$
Then the right hand side of the equation (5.4) is 
$$ h_N(q,z) := {1 \over (q)_N} H_N(q,z).$$
Using the recursion relations (5.10) time after time,
we see that $h_N(q,z)$ satisfies the same recursion relations as (5.8). \qed

\medskip

Next let us evaluate
$$ F_{N,k}(q) = \sum_{\vec{r} \in {\frak R}_N(k)} q^{d(\vec{r})}. $$
Note that the degree of the restricted path $\vec{r}$ could also be calculated as
$$ d(\vec{r}) = \sum_{j=1}^{N-1} j f(r_{N-j}-r_{N+1-j},r_{N-1-j}-r_{N-j}), 
   \eqno(5.11)$$
where 
$$ \eqalign{ f &: \{ +,- \} \times \{ +,- \}  \longrightarrow \{ 0,1 \} \cr 
                   &\cases{f(+,-) = 1 \cr
                      f(+,+) = f(-,-) = f(-,+) = 0.}} \eqno(5.12)$$
By the expression (5.11) of $ d(\vec{r}) $, we can
regard $ d(\vec{r}) $ as a total energy of the size $N+1$ lattice model 
which has the local energy function $f(\;,\;)$. 
This $f(\;,\;)$ is essentially equivalent to the local energy function
of $sl_2$ level $l$ RSOS model. 
Thus we can calculate $ F_{N,k}(q) $ in the same way as
the one dimensional configuration sum of $sl_2$ level $l$ RSOS model [13]. 

\proclaim{Proposition 5.3}.
{$$ F_{N,k}(q) = \sum_{j={- \infty}}^{\infty}
   \Big\{ q^{j(k+1)+{j^2}(l+2)} {N \brack {{-k-2j(l+2)+N} \over 2}}
     -q^{-j(k+1)+{j^2}(l+2)} {N \brack {{k+2-2j(l+2)+N} \over 2}}\Big\}.
      \eqno(5.13) $$}
\par
\hskip-0.7cm{\it Proof.}
Let 
$$\eqalign{
   & g_N^{(\pm,j)}(m,n,n+1) = q^{\pm j(m+1)+j^2 (l+2)} 
     {N \brack {{n \mp(m+1)+1-2j(l+2)+N} \over 2}},\cr
   & g_N^{(\pm,j)}(m,n,n-1) = q^{\pm j(m+1)+j^2 (l+2) + 
     {1 \over 2} \{n \mp (m+1)+1-2j(l+2)+N \}}
     {N \brack {{n \mp(m+1)+1-2j(l+2)+N} \over 2}},\cr
          }\eqno(5.14) $$
and define the function
$$ F_N(m,n,n \pm 1) = \sum_{j={- \infty}}^{\infty}\big\{g_N^{(+,j)}(m,n,n \pm 1)
                        - g_N^{(-,j)}(m,n,n \pm 1) \big\}. \eqno(5.15)$$
Then we can show the following properties:

(1) Initial condition
$$ F_1(m,n,n \pm 1) = g_1^{(+,0)}(m,n,n \pm 1). \eqno(5.16)$$

(2) Recursion relation
$$ \eqalign{
   & F_{N+1}(m,n,n+1) = F_N(m,n-1,n)+F_N(m,n+1,n),\cr
   & F_{N+1}(m,n,n-1) = q^{N+1} F_N(m,n-1,n) + F_N(m,n+1,m),\cr
      } \eqno(5,17) $$

(3) Restriction condition
$$ F_N(m,-1,0) = F_N(m,l+1,l) = 0.\eqno(5.18) $$
From these properties it follows that
$$ F_{N,k}(q) = F_N(a,0,1).\eqno(5.19) $$
Thus we have had the result. \qed

There is also an alternative (fermionic) expression of
$ F_{N,k}(q) $ given by Bouwkneght {\it et al}.[14]

$$ \eqalign{F_{N,k}(q) = \phantom{ q^{-{1 \over 4}k-{1 \over 4}N^2}
    \sum_{m_2,\cdots,m_l} q^{{1 \over 2}
   (N^2+{m_2}^2 +\cdots+{m_l}^2 -N{m_2}-{m_2}{m_3}-\cdots-{m_{l-1}}{m_l})} 
   q^{-{1 \over 4}k-{1 \over 4}N^2}
    \sum_{m_2,\cdots,m_l}ab }\cr
   q^{-{1 \over 4}k-{1 \over 4}N^2}
    \sum_{m_2,\cdots,m_l} q^{{1 \over 2}
   (N^2+{m_2}^2 +\cdots+{m_l}^2 -N{m_2}-{m_2}{m_3}-\cdots-{m_{l-1}}{m_l})}
    \prod_{i=2}^{l} {{1 \over 2}(m_{i-1}+m_{i+1}+ \delta_{i,k+1}) \brack m_i}.\cr}
    $$
where the summation is over all odd non-negative integers for 
$ m_{2a},m_{2a-2},m_{2a-4},\ldots$ and over the even non-negative integers for 
the remaining ones (we set $m_1 = N, m_{l+1}=0 $ ).

\bigskip

\centerline{\bf \S 6. Hidden Yangian symmetry in WZW models}
\bigskip
In this section we point out an intriguing connection between
our spectral decomposition of the path space and the
hidden quantum group symmetry in the WZW conformal
field theory.

The Hilbert space of the $\widehat{sl}_2$ level $l$ WZW model is the direct
sum of the integrable representations $\L(k)$ of $\widehat{sl}_2$ level $l$.
In the level $1$ case the Yangian algebra of $sl_2$ acts on
$\L(k)$ by
$$
Q^a_0=J^a_0,\quad
Q^a_1={1\over2} f^a_{\ bc}\sum_{m>0}J^b_{-m}J^c_m,
\eqno(6.1)
$$
where $J^a_m$ are the Fourier components of the $\widehat{sl}_2$ current and
$f^a_{\ bc}$ is the structure constant of $sl_2$ [1,5].
As a Yangian module, $\L(k)$ decomposes into  irreducible
finite dimensional representations as
$$
\L(k)\cong \bigoplus_{N=k\atop N\equiv k\ {\rm mod}\ 2}^\infty
\bigoplus_{\vec{a}\in Y_N}
W_{N,\vec{a}}.
\eqno(6.2)
$$
Each irreducible component $W_{N,\vec{a}}$ is an $L_0$ eigenspace
with the eigenvalue
$$
{1\over4}(N^2-k^2)+|\vec{a}|+\triangle(k),
\eqno(6.3)
$$ 
and
it has the following $sl_2$-module structure:
$$
W_{N,\vec{a}}\cong \bigotimes_{i=1}^s V_{b_i},
\quad
\beta(\vec{a})=(\beta_1,\dots,\beta_k).
\eqno(6.4)
$$
Comparing (6.2)--(6.4) with Theorem 4.2, we see a remarkable coincidence
of the Yangian decomposition of $\L(k)$ and the spectral decomposition of
the  path space ${\cal P}(k)$.

Let us consider why this happens.
In general   the integrability of field theory is synonymous to
the existence of
 an infinite number of  the local integrals
of motion (IM) commuting with each other.
Let us write the abelian algebra
generated by these IMs as $\cal I$. 
In the spectral decomposition of the path
 space we regard the local energy operators $h_i$
 as the maximal family of the local operators commuting with
 the energy operator $E$ which is also equal to the Virasoro
 energy operator $L_0$. If we identify $h_i$'s
 with the generators
of $\cal I$, then the degeneracy of the spectrum means the presence
of a hidden non-abelian symmetry. In this way we  recover
the decomposition (6.2) through our spectral decomposition.

Let us turn to the higher level case.
We remind that any irreducible finite dimensional
representation of the $sl_2$ Yangian is
isomorphic, as an $sl_2$ module, to a {\it tensor product} of some 
irreducible representations of $sl_2$ [15].
Then, looking at Theorem 4.2, we naturally identify
 each   degeneracy  in our spectral decomposition
with the character of an irreducible Yangian multiplet.
This leads us to the following conjecture:

For $l\geq2$, (for $l=1$, it has been proved in [1])
\proclaim Conjecture. 
1) For each integral representation $\L(k)$ of 
$\widehat{sl}_2$ level $l$, there is
the canonical  action of the Yangian and the
algebra of the local integrals of motion which are commutant
with each other. \hfill\break\indent
2.) The $\L(k)$ decomposes into a direct sum of  irreducible
finite dimensional Yangian modules. The set $
\displaystyle\bigsqcup_{N=k\atop N\equiv k\, {\rm mod}\, 2}^\infty\R_{N}(k)\times Y_N$
parameterizes the Yangian highest weight vectors $v_{\vec{r},\vec{a}}$
in $L(k)$ such  that
the $L_0$ eigenvalue of $v_{\vec{r},\vec{a}}$ is $\Delta(k)+d(\vec{r})
+|\vec{a}|$, and as an $sl_2$-module the Yangian multiplet generated by
$v_{\vec{r},\vec{a}}$ is isomorphic to $\bigotimes_{i=1}^sV_{b_i}$,
where $\beta(\vec{a})=(b_1,\dots,b_s)$.
\par
To some extent the above Yangian module structure   has already appeared
in [14] through
a generalization of the idea of [1].

It has been clarified that the symmetry algebra of
a two dimensional integrable {\it massive} field theory
is also the product of the algebra of the integrals of
motions $\cal I$ and some
``quantum group" symmetry  commuting with each other [16].
Especially, the WZW model allows the integrable massive deformation
which has the Yangian symmetry [17].
It is an important problem to understand how the symmetry algebra of the WZW model
here is related to the one in the deformed model.

We also comment that in the level 1 case there is a simple correspondence
between the  spectrum of the vertex model and the one in
 the Haldane-Shastry model [5]. The natural speculation is that there will be a
 higher spin analog of the Haldane-Shastry model having a similar correspondence 
 to the one in the vertex model here.

To conclude,
we expect that the correspondence between the spectral
decomposition of the solvable lattice model and the one of the corresponding
conformal  field theory  (and its massive deformation)
will be  a universal phenomena.
This interplay 
provides a new and useful way to investigate a hidden quantum symmetry
structure of other conformal field theories as well.
\bigskip

\bigskip
{\bf Acknowledgment}
\bigskip
The authors thank A.~Kuniba, P.~Mathieu, F.A.~Smirnov, and J.~Suzuki for useful discussion.
\bigskip

\bigskip
\centerline{\bf Reference}
\bigskip

\item{[1]}Bernard, D., Pasquier, V., and Serban, D.:  
Spinons in conformal field theory,
Nucl.\ Phys.\ B428, (1994) 612--628

\item{[2]}
Date, E.,   Jimbo,M.,     Kuniba,  A., Miwa, T.,    and Okado, M.:
 Paths, Maya diagrams and representation of $sl(r,C)$,
Advanced Studies in Pure Mathematics 19, (1989) 149--191

\item{[3]} Frenkel, I.\ B.,   and Reshetikhin, N.\ Yu.:   
Quantum affine algebras and holonomic difference equations,
Commun.\ Math.\ Phys.\ 146,  (1992) 1--60

\item{[4]}Tsuchiya, A.,  and  Kanie, Y.: 
Vertex operators in conformal field theory on ${\bf P}^1$
                     and monodoromy representations of braid group,
Adv.\ stud.\ Pure.\ Math.\ 16,  (1988) 297--372

\item{[5]} Haldane, F.~D.~M,  Ha, Z.\ N.\ C.,   Talstra, J.\ C., Bernard, D,
 and Pasquier, V.:
Yangian symmetry of integrable quantum chains
with long range interactions and a new description of states in conformal field theory,
Phys.\ Rev.\ Lett. 69 (1992)
2021--2025

\item{[6]}   Jimbo, M.:
Topics from representation of $U_q(g)$ - an introductory guide to
physicists,
pages 1-61. Nankai Lectures on Mathematical Physics.
World Scientific,  Singapore, 1992

\item{[7]} Baxter, R.\ J.: 
Exactly solvable models in statistical mechanics,
London:  Academic 1982

\item{[8]}  Kang, S-J., Kashiwara, M.,  Misra, K.,  Miwa, T., 
                                       Nakashima, T.,  and Nakayashiki, A.: 
Affine crystals and vertex models,
Int.\ J.\ Mod.\ Phys.\ A  7, Suppl.\ 1A   (1992) 449--484

\item{[9]}   Bouwknegt, P.,  Ludwig, A.,  and Schoutens, K.: 
Spinon basis for $(\hatsl )_k$ integrable highest weight modules
                                             and new character formula,
to appear in Proc.\ of Statistical Mechanics and Quantum Field theory,
USC, May,  (1994) 16--21 (hep-th/9504074)

\item{[10]} Nakayashiki,  A.,  and Yamada,  Y.: 
Crystallizing the spinon basis,
Kyushu Univ.\ preprint  (1995), hep-th/9504052

\item{[11]}
Date,  E.,  Jimbo, M.,     Kuniba, A.,   Miwa, T.,    and Okado, M.: 
One dimensional configuration sums in vertex models and 
affine Lie algebra characters, Lett.~Math.~ Phys.,		
 17 (1989), 69--77

\item{[12]}
Andrews, G.~E.:
The Theory of Partitions,
Addison-Wesley ,(1976)

\item{[13]}
Andrews, G.~E,  Baxter, R.~J,   and  Forrester, P.~J: 
Eight-vertex SOS model and generalized Rogers-Ramanujan type
identities,
J.\ Stat.\ Phys. , 35 (1984) 193--266

\item{[14]}
 Bouwknegt,  P.,   Ludwig, A.,  and Schoutens, K.: 
Spinon bases for higher level SU(2) WZW model,
(hep-th/9412108) 

\item{[15]}
Chari, V.,  and Pressley, A.:  L'Enseignement Math.\ 36
 (1990) 267--302

\item{[16]}
Smirnov, F.A.:  Int.\ J.\ Mod.\ Phys.\ 7A, Suppl.\ 1B
(1992) 813--838, 839--858

\item{[17]}
Bernard, D.:  Commun.\ Math.\ Phys.\ 137 (1991) 191--208
\par

\vskip 2cm
\hskip -0.75cm{\bf Note added.}
After the submission of this paper, 
the authers noticed the  paper,
M.~Idzumi, K.~Iohara,  T.~Jimbo, T.~Miwa, T.~ Nakashima, and  T.~Tokihiro,
{\it Quantum affine symmetry in vertex models} (Int.~J.~Mod.~Phys.~A8(1993),
1479-1511), in which the path space is described using the  combinatorics similar
to that we used for the description of the spectrum of  the path space in section 3.
\bye